\documentclass[12pt,twoside]{article}

\usepackage{amsfonts}
\usepackage{amssymb}
\usepackage{amsmath}
\usepackage{graphicx}

\newcommand{\ignore}[1]{}

\def\@begintheorem#1#2{\par\bgroup{\sc #1\ #2. }\it\ignorespaces}
\def\@opargbegintheorem#1#2#3{\par\bgroup{\sc #1\ #2\ (#3). } \it\ignorespaces}
\def\@endtheorem{\egroup}
\newtheorem{theorem}{Theorem}[section]
\newtheorem{corollary}[theorem]{Corollary}
\newtheorem{lemma}[theorem]{Lemma}
\newtheorem{example}[theorem]{Example}
\newtheorem{proposition}[theorem]{Proposition}
\newtheorem{question}[theorem]{Question}
\newtheorem{definition}[theorem]{Definition}
\newcommand{\bt}[1]{\begin{theorem}\label{#1}}
\newcommand{\bc}[1]{\begin{corollary}\label{#1}}
\newcommand{\bl}[1]{\begin{lemma}\label{#1}}
\newcommand{\be}[1]{\begin{example}\label{#1}}
\newcommand{\bp}[1]{\begin{proposition}\label{#1}}
\newcommand{\bq}[1]{\begin{question}\label{#1}}
\newcommand{\ba}[1]{\begin{algorithm}\rm\label{#1}}
\newcommand{\bd}[1]{\begin{definition}\rm\label{#1}}
\newcommand{\bpr}{\noindent {\em Proof. }}
\newcommand{\et}{\end{theorem}}
\newcommand{\ec}{\end{corollary}}
\newcommand{\el}{\end{lemma}}
\newcommand{\ee}{\end{example}}
\newcommand{\ep}{\end{proposition}}
\newcommand{\eq}{\end{question}}
\newcommand{\ed}{\end{definition}}
\newcommand{\epr}{{\ \vbox{\hrule\hbox{%
\vrule height1.3ex\hskip0.8ex\vrule}\hrule}}\\\par}
\newcommand{\mepr}{{\ \ \ \vbox{\hrule\hbox{%
\vrule height1.3ex\hskip0.8ex\vrule}\hrule}}}

\begin{document}

\title{\bf On Supmodular Matrices}

\author{Shmuel Onn
\thanks{\small Technion - Israel Institute of Technology. Email: onn@technion.ac.il}
}

\date{}

\maketitle

\begin{abstract}
We consider the problem of determining which matrices are permutable to be
supmodular. We show that for small dimensions any matrix is permutable
by a universal permutation or by a pair of permutations, while for higher
dimensions no universal permutation exists. We raise several questions
including to determine the dimensions in which every matrix is permutable.


\end{abstract}

\section{Introduction}

A real $m\times n$ matrix is {\em supmodular} if for all $1\leq i<r\leq m$ and $1\leq j<s\leq n$,
we have $A_{i,j}+A_{r,s}\geq A_{i,s}+A_{r,j}$. Such matrices arise in discrete optimization:
if $A$ is the utility matrix of a transportation problem then an optimal transportation
matrix $X$ is quickly obtained by the greedy algorithm that increases its entries in the order
$X_{1,1},\dots,X_{1,n},\dots,X_{m,1},\dots, X_{m,n}$, each from zero to the maximum possible
value not exceeding prescribed column sums (demands) and row sums (supplies), see \cite{Hof}.
See also \cite{FKP} for related matrix properties under which the greedy algorithm works.

Supmodular matrices are also known as {\em anti-Monge matrices}, see \cite{Fie}.
They are also related to total positivity: a matrix $A$ is supmodular if and only if the matrix
$P:=\exp(A)$ defined by $P_{i,j}:=\exp(A_{i,j})$ is {\em $2$-totally-positive}, namely, all its
minors of order up to $2$ are nonnegative; see \cite{Kar} for this theory and its applications.

Here we are interested in studying which matrices have the following property.
\bd{def}
We say that a real $m\times n$ matrix $A$ is {\em permutable} if its entries
can be permuted in such a way that the permuted matrix is a supmodular matrix.
\ed

If an $m\times n$ matrix $A$ is permutable then so is $A^T$ and so we will assume
$m\leq n$. For an $m\times n$ matrix $\sigma$ whose entries form a permutation of
$1,2,\dots,mn$, let $A^{\sigma}$ be obtained from $A$ by permuting its entries such that
$\sigma_{i,j}<\sigma_{r,s}$ implies $A^{\sigma}_{i,j}\leq A^{\sigma}_{r,s}$.

For instance,
$$
A\ =\
\left(
\begin{array}{ccc}
  1 & 1 & 3 \\
  10 & 3 & 7 \\
  8 & 10 & 6 \\
\end{array}
\right),\quad
\sigma\ =\
\left(
\begin{array}{ccc}
  8 & 7 & 1 \\
  4 & 5 & 3 \\
  2 & 6 & 9 \\
\end{array}
\right),\quad
A^{\sigma}\ =\
\left(
\begin{array}{ccc}
  10 & 8 & 1 \\
  3 & 6 & 3 \\
  1 & 7 & 10 \\
\end{array}
\right).
$$

We prove the following theorem.
\bt{main}
Let $A$ be any $m\times n$ real matrix with $m\leq n$. Then we have:
\begin{enumerate}
\item
If $m=1$ then $A^{\sigma}$ is trivially supmodular for any $\sigma$.
\item
If $m=2$ then $A^{\sigma}$ is supmodular, where
$$\sigma\ =\
\left(
\begin{array}{ccccc}
  n & n-1 & \cdots & 2 & 1 \\
  n+1 & n+2 & \cdots & 2n-1 & 2n \\
\end{array}
\right).
$$
\item
If $m=n=3$ then $A^{\sigma}$ is supmodular, where
$$\sigma\ =\
\left(
\begin{array}{ccc}
  8 & 7 & 1 \\
  4 & 5 & 3 \\
  2 & 6 & 9 \\
\end{array}
\right) .
$$
\item
If $m=3,n=4$ then either $A^{\sigma}$ or $A^{\tau}$ is supmodular, where
$$\sigma\ =\
\left(
\begin{array}{cccc}
  9 &  8 & 7  & 3  \\
  2 &  6 & 5  & 4  \\
  1 & 10 & 11 & 12 \\
\end{array}
\right),\quad
\tau\ =\
\left(
\begin{array}{cccc}
  12 & 3 & 2 & 1  \\
  11 & 7 & 8 & 9  \\
  4  & 5 & 6 & 10 \\
\end{array}
\right).
$$
\end{enumerate}
\et

Theorem \ref{main} asserts that for $m=1,2$ and any $n\geq m$, for $m=n=3$,
and for $m=3,n=4$, any real $m\times n$ matrix $A$ is permutable to a supmodular one.
But moreover, for all cases but the last one, the theorem provides a {\em universal $\sigma$},
that is, one such that $A^{\sigma}$ is supmodular for every $m\times n$
real matrix. Note that the universal permutations are not unique: for $m=1$
all $n!$ permutations are universal, and for $m=2$ and, say, $n=2$, and for $m=n=3$,
the following are universal as well,
$$\sigma\ =\
\left(
\begin{array}{cc}
  4 & 3 \\
  1 & 2 \\
\end{array}
\right),\quad
\sigma\ =\
\left(
\begin{array}{ccc}
  9 & 6 & 2 \\
  3 & 5 & 4 \\
  1 & 7 & 8 \\
\end{array}
\right).
$$
These permutations, as well as these appearing in the theorem, were obtained by using
the notion of goodness of a permutation defined and used in the next section.

Next we show these are the only values of $m,n$ for which a universal $\sigma$ exists.

\bt{universal}
For $m\leq n$, a universal $\sigma$ exists if and only if $m=1,2$ or $m=n=3$.
\et
So already for $m=3,n=4$, none of the $12!=479,001,600$ potential $\sigma$ is universal.

\vskip.2cm
Consider any $m\leq n$. If every  real $m\times n$ matrix $A$ is permutable then
let $p(m,n)\leq(mn)!$ be the smallest positive integer for which there are $\sigma_1,\dots,\sigma_{p(m,n)}$
such that, for every  real $m\times n$ matrix $A$, some $A^{\sigma_i}$ is
supmodular. And if not every $A$ is permutable then let $p(m,n)=\infty$.
Theorem \ref{main} and Theorem \ref{universal} show that
$p(m,n)=1$ if and only if either $m=1,2$ or $m=n=3$, and that $p(3,4)=2$.

\vskip.3cm
Theorem \ref{main} and Theorem \ref{universal} suggest the following question.
\bq{question}\
\begin{enumerate}
\item
What is $p(m,n)$ for all $m\leq n$ and in particular for which $m\leq n$ is it finite?
\item
What is the smallest $m+n$ admitting a non permutable $m\times n$ matrix if any?
\item
What is the complexity of deciding if a given integer matrix is permutable?
\end{enumerate}
\eq

The problem of studying which matrices are permutable is interesting on its own right,
but one possible application is the following. Suppose we need to solve very quickly,
in real time, repeated $m\times n$ transportation problems, with arbitrarily varying
demands and supplies satisfying an upper bound $u$. Suppose we have access to $mn$
transporters, where each transporter $k$ charges $p_k$ per unit flow and can transport
at most $u$ units of flow, so that we cannot simply use only the cheapest.
Our primary objective is to solve the repeated problems very quickly in real time,
and a secondary objective is to solve each with minimum cost. With preprocessing done once
and for all, we try to assign each transporter $k$ to a pair of supplier $i$ and consumer $j$
to transport the flow from $i$ to $j$, so that the resulting utility matrix is supmodular,
and so the repeated problems could be solved very quickly by the greedy algorithm.
This preprocessing is reduced to the problem studied here as follows. We arrange the
negations $-p_k$ of the $mn$ costs arbitrarily in an $m\times n$ matrix $A$, and if $A$
is permutable, search for a permutation $\sigma$ such that $A^\sigma$ is supmodular,
and assign the transporters to pairs $i,j$ according to this permutation. The permuted
matrix $A^\sigma$ is then the utility matrix of all the transportation problems
that we solve (maximizing the utility, so the cost represented by $-A^\sigma$ is minimized),
and all problems can be solve very quickly in real time using the greedy algorithm.

\section{Proofs}

\bl{two_positive}
An $m\times n$ real matrix $A$ is supmodular if and only if we have that, for every $1\leq i<m$
and $1\leq j<n$, the inequality $A_{i,j}+A_{i+1,j+1}\geq A_{i,j+1}+A_{i+1,j}$ holds.
\el
\bpr
Clearly if $A$ is supmodular then the above condition holds. For the converse,
we prove that if the condition holds then, for all $1\leq i<r\leq m$ and
$1\leq j<s\leq n$, we have $A_{i,j}+A_{r,s}\geq A_{i,s}+A_{r,j}$,
by induction on $t=(r-i)+(s-j)$. If $t=2$ this holds by the condition.
Suppose $t>2$ and, say, $s-j>1$. By induction,
\begin{eqnarray*}
(A_{i,j}+A_{r,s})&-&(A_{i,s}+A_{r,j}) \\
&=&\left((A_{i,j}+A_{r,s-1})-(A_{i,s-1}+A_{r,j})\right) \\
&+&\left((A_{i,s-1}+A_{r,s})-(A_{i,s}+A_{r,s-1})\right)
\ \geq\ 0+0\ =\ 0\ .
\mepr
\end{eqnarray*}

We say an $m\times n$ matrix $\sigma$ whose entries are a permutation of $1,2,\dots,mn$
is {\em good on $i,j$} with $1\leq i<m$ and $1\leq j<n$ if
the maximum among $\sigma_{i,j},\sigma_{i,j+1},\sigma_{i+1,j},\sigma_{i+1,j+1}$ is
either $\sigma_{i,j}$ or $\sigma_{i+1,j+1}$ and the
minimum among these entries is either $\sigma_{i,j+1}$ or $\sigma_{i+1,j}$.

An interesting question suggested by a referee is to count and characterize, for every
$m$ and $n$, those permutations which are good for all $1\leq i<m$ and $1\leq j<n$.

\bl{universality_criterion}
Let $\sigma$ be an $m\times n$ matrix whose entries are a permutation of $1,2,\dots,mn$.
If $\sigma$ is good on $i,j$ then
$A^{\sigma}_{i,j}+A^{\sigma}_{i+1,j+1}\geq A^{\sigma}_{i,j+1}+A^{\sigma}_{i+1,j}$
for any real $m\times n$ matrix $A$. If $\sigma$ is not good on $i,j$
then there exists an $A$ with
$A^{\sigma}_{i,j}+A^{\sigma}_{i+1,j+1}<A^{\sigma}_{i,j+1}+A^{\sigma}_{i+1,j}$.
\el
\bpr
First, suppose $\sigma$ is good on $i,j$. Consider any $A$. Then, as claimed,
\begin{eqnarray*}
(A^{\sigma}_{i,j}+A^{\sigma}_{i+1,j+1})&-&(A^{\sigma}_{i,j+1}+A^{\sigma}_{i+1,j}) \\
&=&\left(\max\{A^{\sigma}_{i,j},A^{\sigma}_{i+1,j+1}\}
+\min\{A^{\sigma}_{i,j},A^{\sigma}_{i+1,j+1}\}\right) \\
&-&\left(\max\{A^{\sigma}_{i,j+1},A^{\sigma}_{i+1,j}\}
+\min\{A^{\sigma}_{i,j+1},A^{\sigma}_{i+1,j}\}\right) \\
&=&\left(\max\{A^{\sigma}_{i,j},A^{\sigma}_{i+1,j+1}\}
-\max\{A^{\sigma}_{i,j+1},A^{\sigma}_{i+1,j}\}\right) \\
&+&\left(\min\{A^{\sigma}_{i,j},A^{\sigma}_{i+1,j+1}\}
-\min\{A^{\sigma}_{i,j+1},A^{\sigma}_{i+1,j}\}\right)
\ \geq \ 0+0\ =\ 0 \ .
\end{eqnarray*}
Second, suppose $\sigma$ is not good on some $i,j$, and denote the relevant entries by
$$
\left(
\begin{array}{cc}
  r & p \\
  q & s \\
\end{array}
\right)
\ =\ \left(
\begin{array}{cc}
  \sigma_{i,j} & \sigma_{i,j+1} \\
  \sigma_{i+1,j} & \sigma_{i+1,j+1} \\
\end{array}
\right)\ .
$$
If the minimum among $p,q,r,s$ is $r$ then let $a_1\leq a_2\leq\cdots\leq a_{mn}$ be any
nondecreasing  sequence where $a_r=1$ and $a_p=a_q=a_s=2$. Let $A$ be any matrix
whose entries are $a_1,\dots,a_{mn}$ in any order. Then in $A^{\sigma}$ we have, as claimed,
the inequality
$$(A^{\sigma}_{i,j}+A^{\sigma}_{i+1,j+1})-(A^{\sigma}_{i,j+1}+A^{\sigma}_{i+1,j})
\ =\ (a_r+a_s)-(a_p+a_q)\ =\ -1\ <\ 0\ .$$
A similar argument holds if the minimum is $s$, taking $a_s=1$ and $a_p=a_q=a_r=2$.

If the maximum among $p,q,r,s$ is $p$ then let $a_1\leq a_2\leq\cdots\leq a_{mn}$ be any
nondecreasing  sequence where $a_q=a_r=a_s=1$ and $a_p=2$. Let $A$ be any matrix
whose entries are $a_1,\dots,a_{mn}$ in any order. Then in $A^{\sigma}$ we have, as claimed,
$$(A^{\sigma}_{i,j}+A^{\sigma}_{i+1,j+1})-(A^{\sigma}_{i,j+1}+A^{\sigma}_{i+1,j})
\ =\ (a_r+a_s)-(a_p+a_q)\ =\ -1\ <\ 0\ .$$
A similar argument holds if the maximum is $q$, taking $a_p=a_r=a_s=1$ and $a_q=2$.
\epr

\vskip.3cm\noindent{\em Proof of Theorem \ref{main}.}
Consider any $m\leq n$, any $m\times n$  real matrix $A$, and any
$m\times n$ matrix $\sigma$ whose entries form a permutation of $1,2,\dots,mn$.
Part 1 with $m=1$ holds since any $1\times n$  real matrix is trivially
supmodular.

So assume $m\geq 2$. By Lemma \ref{universality_criterion},
if $\sigma$ is good on $i,j$ for all $1\leq i<m$ and $1\leq j<n$ then
$A^{\sigma}_{i,j}+A^{\sigma}_{i+1,j+1}\geq A^{\sigma}_{i,j+1}+A^{\sigma}_{i+1,j}$ holds
for all such $i,j$ and then $A^{\sigma}$ is supmodular by Lemma \ref{two_positive}.
Part 2 therefore follows since
$$\sigma\ =\
\left(
\begin{array}{ccccc}
  n & n-1 & \cdots & 2 & 1 \\
  n+1 & n+2 & \cdots & 2n-1 & 2n \\
\end{array}
\right)
$$
is good on $1,j$ for all $1\leq j<n$ since the maximum among
$\sigma_{1,j},\sigma_{1,j+1},\sigma_{2,j},\sigma_{2,j+1}$ is $\sigma_{2,j+1}$
and the minimum among these entries is $\sigma_{1,j+1}$.
Part 3 also follows since
$$\sigma\ =\
\left(
\begin{array}{ccc}
  8 & 7 & 1 \\
  4 & 5 & 3 \\
  2 & 6 & 9 \\
\end{array}
\right)
$$
is good on $i=1,2$ and $j=1,2$ as can be verified by direct inspection.

Finally, we prove Part 4. Let $a_1\leq a_2\leq\cdots\leq a_{12}$ be the entries of $A$ arranged
in nondecreasing order. Surely either $a_8+a_5\geq a_7+a_6$ or $a_8+a_5\leq a_7+a_6$ (or both).

First, suppose that $a_8+a_5\geq a_7+a_6$ and consider $A^{\sigma}$ where
$$\sigma\ =\
\left(
\begin{array}{cccc}
  9 &  8 & 7  & 3  \\
  2 &  6 & 5  & 4  \\
  1 & 10 & 11 & 12 \\
\end{array}
\right)\ .
$$
We claim that $A^{\sigma}_{i,j}+A^{\sigma}_{i+1,j+1}\geq A^{\sigma}_{i,j+1}+A^{\sigma}_{i+1,j}$
for $i=1,2$ and $j=1,2,3$ and therefore  $A^{\sigma}$ is supmodular
by Lemma \ref{two_positive}. Indeed, for $i=1,j=2$ this holds since
$$(A^{\sigma}_{1,2}+A^{\sigma}_{2,3})-(A^{\sigma}_{1,3}+A^{\sigma}_{2,2})
\ = \ (a_8+a_5)-(a_7+a_6)\ \geq\ 0\ ,$$
and for every other $i,j$, this follows from Lemma \ref{universality_criterion}
since $\sigma $ is good on $i,j$ as can be verified by inspection. So the claim follows.

Second, suppose that $a_8+a_5\leq a_7+a_6$ and consider $A^{\tau}$ where
$$\tau\ =\
\left(
\begin{array}{cccc}
  12 & 3 & 2 & 1  \\
  11 & 7 & 8 & 9  \\
  4  & 5 & 6 & 10 \\
\end{array}
\right)\ .
$$
We claim that $A^{\tau}_{i,j}+A^{\tau}_{i+1,j+1}\geq A^{\tau}_{i,j+1}+A^{\tau}_{i+1,j}$
for $i=1,2$ and $j=1,2,3$ and therefore  $A^{\sigma}$ is supmodular
by Lemma \ref{two_positive}. Indeed, for $i=2,j=2$ this holds since
$$(A^{\tau}_{2,2}+A^{\tau}_{3,3})-(A^{\tau}_{2,3}+A^{\tau}_{3,2})
\ = \ (a_7+a_6)-(a_8+a_5)\ \geq\ 0\ ,$$
and for every other $i,j$, this follows from Lemma \ref{universality_criterion}
since $\tau$ is good on $i,j$ as can be verified by inspection. So the claim follows.

This completes the proof of Part 4 and the proof of the theorem.
\epr

\vskip.2cm\noindent{\em Proof of Theorem \ref{universal}.}
By Theorem \ref{main} just proved, there exists a universal $\sigma$ for $m=1,2$ and $m=n=3$.
So we need only prove that for all other $m\leq n$ there is no universal $\sigma$.
Suppose for a contradiction that for some $m\geq 3$, $n\geq 4$ there exists a universal $\sigma$,
that is, an $m\times n$ matrix whose entries form a permutation of $1,2,\dots,mn$, such that
$A^{\sigma}$ is supmodular for every  real $m\times n$ matrix $A$.
Let the restriction of $\sigma$ to its top left $3\times 4$ submatrix be
$$
\left(
\begin{array}{cccc}
  a & b & c & d \\
  e & f & g & h \\
  i & j & k & l \\
\end{array}
\right)\ .
$$
Since $\sigma$ is assumed to be universal, by Lemma \ref{universality_criterion}
it must be good on any $2\times 2$ submatrix consisting of consecutive rows and
consecutive columns, that is, among the four entries of such a submatrix,
the maximum must be on the main diagonal and the minimum on the opposite diagonal.

First, suppose $f>g$. Considering $b,c,f,g$, with $f>g$, we obtain $c<g$ and $b>f$.
Considering $c,d,g,h$, with $c<g$, we obtain $h>g$.
Considering $a,b,e,f$, with $b>f$, we obtain $e<f$.
Considering $e,f,i,j$, with $e<f$, we obtain $j>f$.
Considering $f,g,j,k$, with $j>f$, we obtain $k>g$.
Considering $g,h,k,l$, with $k>g$, we obtain $g>h$.
So we obtain the contradiction $g<h<g$.

Second, suppose $f<g$. Considering $f,g,j,k$, with $f<g$, we obtain $j<f$ and $k>g$.
Considering $e,f,i,j$, with $j<f$, we obtain $f<e$.
Considering $a,b,e,f$, with $f<e$, we obtain $b<f$.
Considering $b,c,f,g$, with $b<f$, we obtain $c<g$.
Considering $c,d,g,h$, with $c<g$, we obtain $g<h$.
Considering $g,h,k,l$, with $g<h$, we obtain $k<g$.
So we obtain the contradiction $g<k<g$.
\epr

\section*{Acknowledgments}

Shmuel Onn thanks Steffen Borgwardt for useful related conversations \cite{Bor}.
He was supported by a grant from the Israel Science Foundation and by the Dresner chair.

\end{document}